\documentclass{amsart}
\usepackage{amsfonts}
\usepackage{amsmath}
\usepackage{amssymb}
\usepackage{fullpage}
\usepackage{color}

\newtheorem{theorem}{Theorem}
\newtheorem{corollary}[theorem]{Corollary}
\newtheorem{lemma}[theorem]{Lemma}

\newtheorem{remark}[theorem]{Remark}
\newtheorem{example}[theorem]{Example}

\newcommand{\I}{\mbox{Id}}

\newcommand{\spn}{\mbox{span}}

\DeclareMathOperator{\sgn}{sgn}

\DeclareMathOperator{\var}{var}
\DeclareMathOperator{\Id}{Id}
\DeclareMathOperator{\ad}{ad}
\DeclareMathOperator{\chr}{char}
\DeclareMathOperator{\Der}{Der}
\DeclareMathOperator{\PIexp}{exp}

\begin{document}

\title{Differential identities and polynomial growth of the codimensions}

\author[C. Rizzo]{Carla Rizzo}
\address{Departamento de Matem\'atica, CMUC, Universidade de Coimbra, 3004-504 Coimbra, Portugal}
\email{carlarizzo@mat.uc.pt}

\author[R.B. dos Santos]{Rafael Bezerra dos Santos}
\address{Departamento de Matem\'atica, Instituto de Ci\^encias Exatas, Universidade Federal de Minas Gerais, Avenida Ant\^onio Carlos 6627, 31123-970 Belo Horizonte, Brazil}
\email{rafaelsantos23@ufmg.br}

\author[A.C. Vieira]{Ana Cristina Vieira}
\address{Departamento de Matem\'atica, Instituto de Ci\^encias Exatas, Universidade Federal de Minas Gerais, Avenida Ant\^onio Carlos 6627, 31123-970 Belo Horizonte, Brazil}
\email{anacris@ufmg.br}

\thanks{Corresponding Author: Carla Rizzo.}


\begin{abstract}
	Let $A$ be an associative algebra over a field $F$ of characteristic zero and let $L$ be a Lie algebra over $F$. If $L$ acts on $A$ by derivations, then such an action determines an action of its universal enveloping algebra $U(L)$ and in this case we refer to $A$ as algebra with derivations or $L$-algebra.
	
	Here we give a characterization of the ideal of differential identities of finite dimensional $L$-algebras $A$ in case the corresponding sequence of differential codimensions $c_n^L (A)$, $n\geq 1$, is polynomially bounded. As a consequence, we also characterize $L$-algebras with multiplicities of the differential cocharacter bounded by a constant. 
\end{abstract}

\maketitle

\section{Introduction}
This paper deals with differential identities of algebras over a field $F$ of characteristic zero. More precisely, if $A$ is an associative algebra over $F$ and $L$ is a Lie algebra acting on $A$ by derivations, then this action can be naturally extended to an action of the universal enveloping algebra $U(L)$ of $L$ and in this case we say that $A$ is an algebra with derivations or an $L$-algebra. Then a differential identity of the $L$-algebra $A$ is a polynomial in non-commuting variables $x^d=d(x)$, $d\in U(L)$, vanishing in $A$. Such identities have been studied in later years (see for example \cite{GiambrunoRizzo2018,GordienkoKochetov2014,Kharchenko1979,Rizzo2020ART,Rizzo2021}) and they are a natural generalization of polynomial identities of algebras.

It is well-known that in the ordinary case the polynomial identities satisfied by a given associative algebra $A$ can be measured through its sequence of codimensions $c_n(A)$, $n\geq 1$, i.e., where $c_n(A)$ is the dimension of the space $P_n$ of
multilinear polynomials in $n$ variables modulo the polynomial identities $\I(A)$ of $A$. Such a sequence was introduced by Regev in \cite{Regev1972} and, in characteristic zero, gives an actual quantitative measure of the identities satisfied by a given algebra. The most important feature of the sequence of codimensions proved in \cite{Regev1972} is that if $A$ is an associative algebra satisfying a non trivial polynomial identity (also called PI-algebra), then $c_n(A)$ is exponentially bounded. Later Kemer in \cite{Kemer1979} showed that such codimensions are either polynomially bounded or grow exponentially (no intermediate growth is allowed). In the late nineties Giambruno and Zaicev (see \cite{GiambrunoZaicev1998,GiambrunoZaicev1999}) proved that if $A$ is PI-algebra then the limit $\PIexp(A):=\lim_{n\to \infty}\sqrt[n]{c_n(A)}$ exists and is always a non-negative integer called the exponent of $A$.

One of the main advantages of the exponent is to have an integral scale allowing us classify the PI-algebras according to their exponent. Much effort has been put into the study PI-algebras with “slow” codimension growth. It is well known that $\PIexp(A) \leq 1$ if and only if $c_n(A)$ is polynomially bounded. Various descriptions of such algebras were given (see for example \cite{GiambrunoZaicev2001, Kemer1978, Kemer1979}). Similar results were also proved in the setting of varieties of graded algebras and algebras with involution (see for example \cite{GiambMish, GiambMishZai,KoshlukovLaMattina2015, Valenti2011}).

Inspired by the above results it is natural to expect that  similar conclusions hold for varieties of $L$-algebras. In fact, in analogy with the ordinary case, one defines the
sequence of differential codimensions $c_n^L(A)$, $n\geq 1$, of an $L$-algebra $A$ and the growth of the $L$-variety $\mathcal{V}= \var ^{L}(A)$, i.e., variety of algebras with derivations, is
the growth of the sequence $c_{n}^{L}(\mathcal{V})=c_{n}^{L}(A)$,
$n=1,2,\dots$. In case $A$ is a finite dimensional $L$-algebra, Gordienko in \cite{Gordienko2013JA} proved that the limit $\PIexp^L(A):=\lim_{n\to \infty}\sqrt[n]{c_n^L(A)}$ exists and is a non-negative integer called $L$-exponent of $A$. As a consequence, it follows that the differential codimensions of a finite dimensional algebra are either polynomially bounded or grow exponentially.

Our purpose here is to characterize $L$-varieties $\mathcal{V}$ having polynomial growth and we reach our goal in the setting of varieties generated
by finite dimensional $L$-algebras $A$. In this frame we prove that $\PIexp^L(A)\leq 1$ if and only if $\PIexp (A)\leq 1$. As a consequence, we get that $c_n^L(A)$ is polynomially bounded if and only if $c_n(A)$ is polynomially bounded. Notice that this result was known before only in case $L$ is a finite dimensional semisimple Lie algebra (see \cite[Theorem 15]{GordienkoKochetov2014}).

Similarly to the ordinary case, another two useful invariants can be attached to an algebra with derivations $A$: the sequence of differential cocharacter $\chi_n^L(A)$, $n \geq 1$, where $\chi_n^L(A)$ is the character of the $S_n$-module of
multilinear differential polynomials in $n$ variables modulo the differential identities $\I^L(A)$ of $A$, and the differential colength sequence $l_n^L(A)$, $n\geq 1$, where $l_n^L(A)$ is the sum of the corresponding multiplicities of $\chi_n^L(A)$.

It is well-known that, in case $A$ is a finite dimensional $L$-algebra, the multiplicities of the differential cocharacter are polynomially bounded (see \cite{Gordienko2013JA}). Thus it seems interesting also to characterize the differential cocharacter sequence when stronger conditions hold for the multiplicities. In this perspective, motivated by the results for ordinary algebras \cite{MischenkoRegevZaicev1999}, for graded algebras \cite{CirritoGiambruno2013,Otera2005} and for algebras with involution \cite{Rizzo2020JA,Vieira2015}, we characterize the differential identities when the corresponding multiplicities are bounded by a constant. In particular we prove that the multiplicities in $\chi_n^L(A)$ are bounded by a constant if and only if differential codimensions of $A$ grow polynomially, and, consequently, we get another characterization of $L$-varieties of polynomial growth. Also as a direct consequence of this results we have that $c_n^L(A)$ is polynomially bounded if and only if $l_n^L(A)$ is bounded by a constant. 

We give also two others characterizations of $L$-varieties $\mathcal{V}$ of polynomial growth: the first one in terms of the structure of an algebra generating $\mathcal{V}$ and the second one in terms of the shape of the diagrams of the irreducible $S_n$-characters appearing with non-zero multiplicity in the $n$th differential cocharacter of $\mathcal{V}$.

\section{Preliminaries}
Throughout this paper $F$ will denote a field of characteristic zero and $L$ a Lie algebra over $F$. Let $A$ be an associative algebra over $F$. Recall that a derivation of $A$ is a linear map $\delta:A\to A$ such that
$$\delta(ab)=\delta(a) b+a \delta(b), \mbox{ for all } a,b\in A.$$
In particular, an inner derivation induced by $ a\in A $ is the derivation $\ad_a:A\to A$ of $A$ defined by $\ad_a (b)=[a,b]=ab-ba$, for all $b\in A$. The set of all derivations of $A$ is a Lie algebra denoted by $\Der(A)$, and the set $\ad (A)$ of all inner derivations of $A$ is a Lie subalgebra of $\Der(A)$.

If $L$ acts on $A$ by derivations, then by the Poincaré-Birkhoff-Witt Theorem, the $L$-action on $A$ can be naturally extended to an $U(L)$-action, where $U(L)$ is the universal enveloping algebra of $L$. In this way $A$ becomes a left $U(L)$-module and we call it algebra with derivations or $L$-algebra.

Given a basis $\mathcal{B}=\{h_{i}: i\in I \}$ of $U(L)$, we let $F\langle X|L\rangle$ be the free associative algebra over $F$ with free formal generators $x_{j}^{h_{i}}$, $i\in I$, $j\in\mathbb{N}$. 
For all $h=\sum_{i\in I}\alpha_i h_i\in U(L)$, where only a finite number of $\alpha_i\in F$ are non-zero, we set $x^h:= \sum_{i\in I}\alpha_i x^{h_i}$. We let $U(L)$ act on $F\langle X| L\rangle$  by setting
$$\gamma(x_{j_{1}}^{h_{i_{1}}}x_{j_{2}}^{h_{i_{2}}}\dots x_{j_{n}}^{h_{i_{n}}})=x_{j_{1}}^{\gamma h_{i_{1}}}x_{j_{2}}^{h_{i_{2}}}\dots x_{j_{n}}^{h_{i_{n}}}+\dots+x_{j_{1}}^{h_{i_{1}}}x_{j_{2}}^{h_{i_{2}}}\dots x_{j_{n}}^{\gamma h_{i_{n}}},$$
where $\gamma\in L$ and $x_{j_{1}}^{h_{i_{1}}}x_{j_{2}}^{h_{i_{2}}}\dots x_{j_{n}}^{h_{i_{n}}}\in F\langle X|L\rangle$. In this way $F\langle X| L\rangle$ has a structure of $L$-algebra. We write $x_{i}:=x_{i}^{1}$, $1\in U(L)$, and we set $X=\{x_{1},x_{2},\dots\}$.
Then the algebra $F\langle X|L\rangle$ is called the free associative algebra with derivations on the countable set $X$ over $F$ and its elements are called differential polynomials.

Let now $A$ be an $L$-algebra. A polynomial $f(x_{1},\dots,x_{n})\in F\langle X|L\rangle$ is a differential identity of $A$, or an $L$-identity of $A$, if $f(a_{1},\dots,a_{n})=0$ for all $a_{i}\in A$, and, in this case, we write $f\equiv 0$. We denote by
$$
\Id^{L}(A)=\{f\in F\langle X|L\rangle : f\equiv 0 \mbox{ on } A\},
$$
the $T_{L}$-ideal of differential identities of $A$, i.e., it is an ideal of $F\langle X|L\rangle$ invariant under the $U(L)$-action. In characteristic zero $\Id^L(A)$ is completely determined by its multilinear polynomials and for every $n\geq 1$ we denote by
$$
P_{n}^{L}=\spn\{x_{\sigma(1)}^{h_{i_1}}\dots x_{\sigma(n)}^{h_{i_n}}: \sigma\in S_{n},h_{i}\in \mathcal{B} \}
$$
the space of multilinear differential polynomials of degree $n$. Notice that in case $A$ is a finite dimensional $L$-algebra, $U(L)$ acts on $A$ as a suitable finite dimensional subalgbera of the endomorphism algebra of $A$, then we may assume that $
P_{n}^{L}$ is finite dimensional and similarly to the ordinary case we can define the following invariants.

Let us consider the space
$$P_n^L(A)=\dfrac{P_{n}^{L}}{P_{n}^{L}\cap \Id^{L}(A)}, \quad n\geq 1.$$
Then its dimension $c_n^L(A)=\dim_F P_n^L(A)$
is called the $n$th differential codimension of $A$. Moreover recall that the symmetric group $S_n$ acts on the left on the space $P_{n}^{L}$ as follows: for $\sigma\in S_{n}$, $\sigma(x_{i}^{h})=x_{\sigma(i)}^{h}$. Since $P_{n}^{L}\cap \Id^{L}(A)$ is stable under this $S_{n}$-action, then $P_n^L(A)$
is a left $S_{n}$-module and its character, denoted by $\chi_{n}^{L}(A)$, is called the $n$th differential cocharacter of $A$. Since $F$ is of characteristic zero, we can write
$$\chi_n^{L}(A)=\sum_{\lambda\vdash n} m_{\lambda}\chi_{\lambda},$$
where $\lambda$ is a partition of $n$, $\chi_{\lambda}$ is the irreducible $S_{n}$-character associated to $\lambda$ and $m_{\lambda}\geq 0$ is the corresponding multiplicity.

Another numerical sequence that can be attached to a $L$-algebra $A$ is the sequence of differential colengths. If $\chi_n^{L}(A)=\sum_{\lambda\vdash n} m_{\lambda}\chi_{\lambda}$ is the $n$th differential cocharacter of $A$, then the $n$th differential colength of $A$ is defined as
$$
l_n^L(A)= \sum_{\lambda\vdash n} m_{\lambda}.
$$

Notice that the theory of differential identities generalizes the ordinary theory of polynomial identities. In fact, any algebra $A$ can be regarded as $L$-algebra by letting $L$ act on $A$ trivially,
i.e., $L$ acts on $A$ as the trivial Lie algebra and $U(L)\cong F$.
Moreover, since $U(L)$ is an algebra with unit, we can identify in a natural way $P_n$ with a subspace of $P_n^L$. Hence $P_n\subseteq P^L_n$ and $P_n\cap \Id(A)=P_n\cap \Id^L(A)$. If we consider the $S_n$-modulo $P_n(A)=\dfrac{P_{n}}{P_{n}\cap \Id(A)}$, $n\geq 1$, then its dimension $c_n(A)=\dim_F  P_n(A)$ is the $n$th (ordinary) codimension of $A$ and its character $\chi(A)$ is the $n$th (ordinary) cocharacter of $A$.  Since $ \chr F=0$, by complete reducibility we can write 	
$$
\chi_n(A)=\sum_{\lambda\vdash n} \bar{m}_\lambda \chi_\lambda,
$$
where $\lambda$ is a partition of $n$, $\chi_{\lambda}$ is the irreducible $S_{n}$-character associated to $\lambda$ and $\bar{m}_{\lambda}\geq 0$ is the corresponding multiplicity.
As a consequence we have the following relations.
\begin{remark}\label{Rmk Codimensions and multiplicity}
	For all $n\geq 1$,
	\begin{itemize}
		\item[1.]  $c_n(A)\leq c_n^L(A)$;
		\vspace{1mm}
		\item[2.]  $\bar{m}_\lambda\leq m_\lambda$, for any $\lambda\vdash n$.
	\end{itemize}
\end{remark}

Recall that if $A$ is an $L$-algebra, then the variety of algebras with derivations generated by $A$ is denoted by $\var ^{L}(A)$ and is called $L$-variety. The growth of $\mathcal{V}= \var ^{L}(A)$ is the growth of the sequence $c_{n}^{L}(\mathcal{V})=c_{n}^{L}(A)$, $n\geq 1$, then we say that the $L$-variety $\mathcal{V}$ has polynomial growth if $c_{n}^{L}(\mathcal{V})$ is polynomially bounded. In what follow we shall characterize $L$-variety of polynomial growth.

We conclude this section by recalling a basic result concerning the sequence of cocharacters which can be easily proved.

\begin{remark} \label{Remark cocharacter and colegth}
	Let $A$ and $B$ be two $L$-algebras such that
	$$\chi_n^{L}(A)=\sum_{\lambda\vdash n} m_{\lambda}\chi_{\lambda} \quad \mbox{and} \quad \chi_n^{L}(B)=\sum_{\lambda\vdash n} m_{\lambda}' \chi_{\lambda}.
	$$
	Then the direct sum $A\oplus B$ is also an $L$-algebra where the $L$-action is induced by the $L$-action by derivations defined on $A$ and $B$. Moreover, if
		$$\chi_n^{L}(A\oplus B)=\sum_{\lambda\vdash n} \tilde{m}_{\lambda}\chi_{\lambda}
		$$
		is the decomposition of the $n$th differential cocharacter of $A\oplus B$, then $\tilde{m}_\lambda \leq m_\lambda + m_\lambda '$, for all $\lambda\vdash n$.

\end{remark}

\section{Finite dimensional $L$-algebras and varieties of polynomial growth}

In this section we shall characterize finite dimensional algebras with derivations generating varieties of polynomial growth.

We start by recalling some results on the structure of finite dimensional algebras with derivations.

Let $L$ be a Lie algebra over $F$ and $A$ an $L$-algebra over $F$. An ideal (subalgebra) $I$ of $A$ is an $L$-ideal (subalgebra) if it is an ideal (subalgebra) such that $I^L\subseteq I$, where $I^L$ denotes the set of all $h(a)$, for all $a\in I$ and $h\in U(L)$. The algebra $A$ is $L$-simple if $A^2\neq \{0\}$ and $A$ has no non-trivial $L$-ideals.

Let $A$ be a finite dimensional $L$-algebra over $F$. By Wedderburn-Malcev Theorem for associative algebras (see \cite[Theorem 3.4.3]{GiambrunoZaicev2005book}), we can write $A$ as a direct sum of vector spaces
$$A=A_{ss}+J,$$
where $A_{ss}$ is a maximal semisimple subalgebra of $A$ and $J=J(A)$ is the Jacobson radical of $A$. Notice that $J$ is always an $L$-ideal of $A$ (see \cite[Theorem 4.2]{Hoch}), but it may not exist an $L$-invariant Wedderburn-Malcev decomposition, i.e., it may happen that $A_{ss}^L\nsubseteq A_{ss}$, for every maximal semisimple subalgebra $A_{ss}$ of $A$. However, we remark that the Wedderburn-Malcev decomposition always exists in case $L$ is a semisimple Lie algebra (see \cite[Theorem 4]{GordienkoKochetov2014}). In what follows we give an example of an $L$-algebra that has no $L$-invariant Wedderburn-Malcev decomposition.

\begin{example}
	Let $UT_2^\delta$ be the $L$-algebra of $2 \times 2$ upper triangular matrices where $L$ acts on it as the 1-dimensional Lie algebra spanned by the inner derivation $\delta=\ad_{e_{12}}$. Suppose that $UT_2 ^\delta =A_{ss}+J$ for some maximal semisimple subalgebra $A_{ss}$  of $UT_2^\delta$ such that $A_{ss}^L\subseteq A_{ss}$. Since $\delta=\ad_{e_{12}}$, $J=\spn_{F}\{e_{12}\}$ and $A_{ss}^L\subseteq A_{ss}$, it follows that $[A_{ss},J]\subseteq A_{ss}$. On the other hand , since $J$ is an ideal of $UT_2^\delta$, $[A_{ss},J]\subseteq J$. Thus it follows that $[A_{ss},J]\subseteq A_{ss}\cap J =\{0\}$. But since $J=\spn_{F}\{e_{12}\}$, we have that $[J,J]=\{0\}$. This says that the center of $UT_2^\delta$ contains $J$, that is no true. Therefore $A_{ss}^L\nsubseteq A_{ss}$, for all maximal semisimple subalgebra $A_{ss}$ of $UT_2^\delta$. Thus $UT_2^\delta$ has no $L$-invariant Wedderburn-Malcev decomposition.
\end{example}

In \cite{Gordienko2013JA}, Gordienko proved that if $A$ is a finite dimensional $L$-algebra, then the sequence of differential codimensions $c_n^{L}(A)$ is exponentially bounded. Moreover, the author proved that the limit 
$\displaystyle \lim_{n\to \infty}\sqrt[n]{c_n^L(A)}
$
exists and is a non-negative integer. In this case, this limit is called the $L$-exponent of $A$ and is denoted by $\exp^L(A)$. In particular, we have the following.

\begin{theorem}\label{exp}\cite[Theorems 1 and 3]{Gordienko2013JA} Let $A$ be a finite dimensional algebra over a field of characteristic zero. If $L$ is a Lie algebra acting on $A$ by derivations, then there exist constants $C_1,C_2,r_1,r_2, C_1>0,$ and a positive integer $d$ such that
	$$C_1n^{r_1}d^n\leq c_n^L(A)\leq C_2n^{r_2}d^n, \mbox{ for all } n\in \mathbb{N}.$$ 
	Hence, $\exp^L(A)=d$. Moreover, If $J=J(A)$ is the Jacobson radical of $A$ and $A/J=\overline{A}_1\oplus \cdots \oplus \overline{A}_m$, then
	$$d=\max\{\dim(\overline{A}_{i_1}\oplus \overline{A}_{i_2}\oplus\cdots\oplus \overline{A}_{i_k}): A_{i_1}^L A^+  A_{i_2}^L A^+ \cdots A^+  A_{i_k}^L\neq \{0\}\},$$ 
	where $i_r\neq i_s, 1\leq r,s\leq n$, $A^+=A+ F\cdot 1$ and $A_i$ is a subalgebra of $A$ (not necessary $L$-invariant) such that $\pi(A_i) =\overline{A}_i$, for all $1\leq i \leq m$, where $\pi: A\to A/J$ is the natural projection.
	
\end{theorem}

As a consequence we have the following.

\begin{corollary}
	If $A$ is a finite dimensional $L$-algebra, the sequence $c_n^L(A)$, $n\geq 1$, either is polynomially bounded or growth exponentially.
\end{corollary}

\begin{lemma}
	\label{A_1JA_2}
	Let $A$ be a finite dimensional $L$-algebra over an algebraically closed field $F$ of characteristic zero such that, as ordinary algebra, $A= A_1 \oplus \dots \oplus A_m +J$ with $A_1 \cong \dots \cong A_m \cong F$.  If there exist $1 \leq i,k\leq m$, $i \neq k$, such that $A_i ^L A^+ A_k ^L \neq \{0\}$, then  $A_i J A_k \neq \{ 0\}$.
\end{lemma} 
\begin{proof}
	We assume, as we may, that $i=1$ and $k=2$.
	Let $a\in A^+$, $e_i\in A_i$ with $e_i^2=e_i$, $i=1,2$, such that $h_1(e_1) a h_2(e_2)\neq 0$, for some $h_1,h_2\in U(L)$. If $h_1, h_2\in \spn_F \{1_{U(L)}\}$, then we are done. So let us suppose that $h_1\notin \spn_F \{1_{U(L)}\}$ and $h_2\in \spn_F \{1_{U(L)}\}$, the other cases will follow analogously. Without loss generality we may suppose that $h_1=\gamma_1 \dots \gamma_r$, $\gamma_i\in L$, $i=1,\dots,r$, $r\geq 1$, and $h_2=1_{U(L)}$.  We proceed by induction on $r$.
	
	If $r=1$, then $h_1=\gamma_1\in L$, and since $e_1^2=e_1$, we have that $h_1(e_1)a e_2= e_1 h_1(e_1)a e_2+ h_1(e_1)e_1 a e_2\neq 0$. Thus if $e_1h_1(e_1)a e_2\neq 0$,  we are done since $h_1(e_1)a\in J$. If $h_1(e_1) e_1 a e_2 \neq 0$, then $a\in J$ and we are done. So let $r>1$. Hence by definition of derivation and the idempotence of $e_1$ we have that
$$
	h_1(e_1)=h_1(e_1)e_1+ e_1 h_1(e_1)+ \sum_{I,K} h_I(e_1) h_K(e_1),
$$
	where $I=\{i_1, \dots, i_p\}$ and $K=\{k_1, \dots, k_t\}$ are two disjoint subsets of $\{1, \dots, r\}$ such that $i_1<\dots < i_p$, $p<r$, and $k_1<\dots< k_t$, $t<r$, respectively, $h_I= \gamma_{i_1}\cdots \gamma_{i_p}$ and $h_K=\gamma_{k_1}\cdots \gamma_{k_t}$. Since $h_1(e_1)a e_2\neq 0$, then
	$$
	h_1(e_1)e_1a e_2+ e_1 h_1(e_1)a e_2+ \sum_{I,K} h_I(e_1) h_K(e_1)a e_2\neq 0.
	$$
	
	If $e_1h_1(e_1)a e_2\neq 0$ or $h_1(e_1)e_1 a e_2\neq 0$, then we are done. So let us assume that there exist $I$ and $K$ such that $h_I(e_1) h_K(e_2)a e_2\neq 0$. Then it follows that $h_K(e_2)a e_2\neq 0$ and by inductive hypothesis the proof is complete.
	
\end{proof}

If we denote by $\PIexp(A)=\displaystyle \lim_{n\to \infty}\sqrt[n]{c_n(A)}$ the ordinary exponent of $A$, then
$$
\PIexp(A)=\max \{\dim(A_{i_1}\oplus \dots \oplus A_{i_k}) : A_{i_1}JA_{i_2}J\dots JA_{i_k}\neq \{0\}\},
$$
where $J=J(A)$ is the Jacobson radical of $A$ and $A= A_1 \oplus \dots \oplus A_m +J$ is the Wedderburn-Malcev decomposition of $A$ as ordinary algebra (see \cite[Chapter 6]{GiambrunoZaicev2005book}). Thus as an immediate consequence of Lemma \ref{A_1JA_2} and Theorem \ref{exp} we have the following.

\begin{theorem}
	\label{Thm: Characterization in terms of ordinary exp}
	Let $L$ be a Lie algebra over a field $F$ of characteristic zero and $A$ be a finite dimensional $L$-algebra over $F$.  Then the following conditions are equivalent:
	\begin{itemize}
		\item[1.] $c_n^L(A)$ is polynomially bounded;
		\vspace{1mm}
		\item[2.] $\PIexp^L(A)\leq 1$;
		\vspace{1mm}
		\item[3.] $c_n(A)$ is polynomially bounded;
		\vspace{1mm}
		\item[4.] $\PIexp(A)\leq 1$.
	\end{itemize}
\end{theorem}

As in the ordinary case, we have the following remark (see \cite[Lemma 7.2.1]{GiambrunoZaicev2005book}).

\begin{remark}
	If $A$ and $B$ are $L$-algebras, then $A\oplus B$ has an induced structure of $L$-algebra and $c_n^L(A\oplus B)\leq c_n^{L}(A)+c_n^{L}(B)$. As a consequence, $\exp^L(A\oplus B)=\max\{\exp^L(A),\exp^{L}(B)\}$.
\end{remark}

 Recall that if $A$ and $B$ are two $L$-algebras, then we say that $A$ is $T_L$-equivalent to $B$, and we write $A\sim_{T_L}B,$ if $\Id^L(A)=\Id^L(B)$. Notice that given an $L$-algebra $A$, $A$ is $T_L$-equivalent to $B$ if and only if $\var^L (A)=\var^L (B)$.

\begin{lemma}
	\label{teclemma}
	Let $F$ be a field of characteristic zero, $\bar{F}$ the algebraic closure of $F$ and $A$ a finite dimensional $L$-algebra over $\bar{F}$, where $L$ is a Lie algebra over $\bar{F}$ acting on $A$ by derivations. Suppose that $\dim_{\bar{F}} A/J(A)\leq 1$. Then $A\sim_{T_L} B$ for some finite dimensional $L$-algebra $B$ over $F$ with $\dim_F B/J(B)\leq 1$.
\end{lemma}

\begin{proof}
	Since $\dim_{\bar{F}} A/J(A)\leq 1$, it follows that either $A\cong \bar{F}+J(A)$ or $A = J(A)$ is a nilpotent algebra. Now we take an arbitrary basis $\{v_1,\dots, v_p\}$ of $J(A)$ over $\bar{F}$ and we let $B$ be the $L$-algebra over $F$ generated by $\mathcal{B}=\{1_{\bar{F}}, v_1, \dots, v_p\}$ or $\mathcal{B}=\{ v_1, \dots, v_p\}$ according as $A\cong \bar{F}+J(A)$ or $A = J(A)$, respectively.
	
	Since $A$ is finite dimensional over $\bar{F}$ and $J(A)$ is a nilpotent $L$-ideal of $A$, $B$ is finite dimensional over $F$. Therefore $B$ is a finite dimensional $L$-algebra and $\dim_F B/J(B)=\dim_{\bar{F}} A/J(A)\leq 1$. Now notice that, as $F$-algebras, $\Id^L (A)\subseteq \Id^L(B)$. On the other hand, if $f$ is a multilinear differential identity of $B$ then $f$ vanishes on $\mathcal{B}$. But $\mathcal{B}$ is a basis of $A$ over $\bar{F}$. Hence $ \Id^L(B) \subseteq \Id^L (A)$ and $A \sim_{T_L} B$. 
\end{proof}

Next theorem gives a characterization of $L$-varieties of polynomial growth in terms of the structure of the generating algebra.

\begin{theorem}
	\label{Thm: Characterization in terms of equvalent algebra}
	Let $L$ be a Lie algebra over a field $F$ of characteristic zero and $A$ be a finite dimensional $L$-algebra over $F$. Then $c_n^L(A)$, $n\geq 1$, is polynomially bounded if and only if $A\sim_{T_L} B_1 \oplus \dots \oplus B_m$, where $B_1,\dots, B_m$ are finite dimensional $L$-algebras over $F$ such that $\dim B_i/J(B_i)\leq 1$ for all $1\leq i \leq  m$.
\end{theorem}

\begin{proof}
	Suppose first that $A\sim_{T_L}B$ where $B=B_1 \oplus \dots \oplus B_m,$ with $B_1,\dots, B_m$ finite dimensional $L$-algebras over $F$ such that $\dim B_i/J(B_i)\leq 1$, for all $1\leq i \leq  m$. Then, by Theorem \ref{exp}, $c_n ^L (B_i)$ is polynomially bounded, for all $1\leq i \leq  m$, and $c_n ^L (A)=c_n ^L (B)\leq c_n ^L (B_1)+ \cdots + c_n ^L (B_m)$. Thus $c_n ^L (A)$ is polynomially bounded.
	
	Conversely, suppose that $c_n^L(A)$ is polynomially bounded. Assume first that $F$ is algebraically closed. Let $A=A_{ss}+J$ where $A_{ss}$ is a semisimple subalgebra and $J=J(A)$ is the Jacobson radical of $A$. 
	By Theorem \ref{exp}, it follows that $A_{ss}=A_1 \oplus \dots \oplus A_l$ with $A_1 \cong \cdots \cong A_l \cong F$ and $A_i ^L A^+ A_k ^L=\{0\}$, for all $1 \leq  i, k \leq l$, $i\neq k$.
	
	Set $B_1=A_1+ J, \dots, B_l= A_l+J$. Since $A_i^L \subseteq A_i +J$ for all $1\leq i \leq  l$, and $J$ is an $L$-ideal of $A$, $B_i$ is an  $L$-subalgebra of $A$, for all $1\leq i \leq  l$. We claim that $$\Id^{L}(A)=\Id^{L}(B_1) \cap \dots \cap \Id^L( B_l) \cap \Id^L( J).$$

	Clearly $\I^L(A)\subseteq \Id^{L}(B_1) \cap \dots \cap \Id^L( B_l) \cap \Id^L( J)$.
	Now let $f\in \Id^{L}(B_1) \cap \dots \cap \Id^L( B_l) \cap \Id^L( J)$ and suppose that $f$ is not a differential identity of $A$. We may clearly assume that $f$ is multilinear. Moreover, by choosing a basis of $A$ as the union of a basis of $A_{ss}$ and a basis of $J$, it is enough to evaluate $f$ on this basis. Let $u_1,\dots,u_t$ be elements of this basis such that $f(u_1,\dots, u_t)\neq 0$. Since $f\in \Id^L (J),$ at least one element, say $u_s$, does not belong to $J$. Then $u_s\in B_r$, for some $r$. Recalling that $A^L_iA^L_k\subseteq A_i ^L A^+ A_k ^L=\{0\}$, for all $i\neq k$, we must have that $u_1,\dots, u_t\in A_r \cup J$. Thus $u_1,\dots,u_t\in A_r + J=B_r$ and this contradicts the fact that $f$ is a differential identity of $B_r$. This prove the claim. The proof is completed by noticing that $\Id^{L}(B_1 \oplus \dots \oplus B_l \oplus J)=\Id^{L}(B_1) \cap \dots \cap \Id^L( B_l) \cap \Id^L( J)$ and $\dim B_i/J(B_i)=1$, for all $1\leq i \leq  l$.
	
	In case $F$ is arbitrary, we consider the algebra $\bar{A}=A \otimes_F \bar{F}$, where $\bar{F}$ is the algebraic closure of $F$ and $\bar{A}=A \otimes_F \bar{F}$ is endowed with the induced $L$-action $(a\otimes\alpha)^\gamma=a^\gamma \otimes \alpha$, for $\gamma\in L$, $a\in A$ and $\alpha \in \bar{F}$. Clearly, over $F$, $\var^L (A)=\var^L (\bar{A})$. Moreover, the differential codimensions of $A$ over $F$ coincide with the differential codimensions of $\bar{A}$ over $\bar{F}$. Thus, by hypothesis, it follows that the differential codimensions of $\bar{A}$ are polynomially bounded. But then, by the first part of the proof, $\bar{A}\sim_{T_L} B_1 \oplus\dots \oplus B_m$ where $B_1,\dots, B_m$ are finite dimensional $L$-algebras over $\bar{F}$ such that $\dim_{\bar{F}} B_i/J(B_i)\leq 1$, for all $1\leq i \leq  m$. By Lemma \ref{teclemma} there exist finite dimensional $L$-algebras $C_1,\dots, C_m$ over $F$ such that, for all $i$, $C_i \sim_{T_L} B_i$ and $\dim_F C_i / J(C_i) \leq 1$. It follows that $\Id^L(A)=\Id^L(\bar{A})=\Id^L (B_1 \oplus\dots \oplus B_m)=\Id^L (C_1 \oplus\dots \oplus C_m)$ and we are done.
\end{proof}

\section{Differential cocharacter of varieties of polynomial growth}

In this section we give other characterizations of $L$-varieties $\mathcal{V}$ of polynomial growth through the behaviour of their sequences of cocharacters.

\begin{theorem}
	\label{Thm: Characterization in terms of cocharacter}
	Let $L$ be a Lie algebra over a field $F$ of characteristic zero and let $A$ be a finite dimensional $L$-algebra over $F$. Then $c_n^L(A)$, $n\geq 1$, is polynomially bounded if and only if there exists a constant $q$ such that
	$$\chi_n^L(A)=\sum_{\substack{\lambda\vdash n \\ |\lambda|-\lambda_1< q}} m_\lambda \chi_\lambda$$
	and $J(A)^q=\{0\}$.
\end{theorem}

\begin{proof}
	Notice that the decomposition of $\chi_n^L(A)$ into irreducible characters does not change under extensions of the base field. This fact can be proved following word by word the proof for the ordinary case (see for example \cite[Theorem 4.1.9]{GiambrunoZaicev2005book}). Also if $\bar{F}$ is the algebraic closure of $F$ and $J(A)^q=\{0\}$, then $J(A\otimes_F \bar{F})^q=\{0\}$. Therefore we may assume, without loss of generality, that $F$ is an algebraically closed field.
	
	Suppose $c_n^L(A)$, $n\geq 1$, is polynomially bounded and let $\lambda$ be a partition of $n$ such that $|\lambda|-\lambda_1 \geq q$ and $m_\lambda \neq 0$. Then there exist $f\in P_n ^L$ and a tableau $T_\lambda$ such that $e_{T_\lambda} f\notin \Id^L (A)$. Let $\lambda' =(\lambda_1', \dots,  \lambda_t ')$ be the conjugate partition of $\lambda$. Then $e_{T_\lambda} f$ is a linear combination of polynomials each alternating on $t$ disjoint sets of $\lambda_1 '  , \dots, \lambda_t '$ variables, respectively. We shall reach a contradiction by proving that these polynomials $g$ vanish in $A$.

	Let $A=A_1 \oplus \dots \oplus A_m +J$, where $A_1, \dots, A_m$ are simple algebras and $J=J(A)$ is the Jacobson radical, then by Theorem \ref{exp}, $\dim A_i =1$ and $A_i ^L A^+ A_k ^L =\{0\}$ for all $1 \leq i,k\leq m$, $i\neq k$.
	 In order to get a non-zero value of $g$ we must replace at most one variable with elements of a single component, say, $A_i$, and the others variables with elements of $J$. Since $\dim A_i=1$, we can substitute at most one element of $A_i$ in each alternating set. Thus we can substitute at most $\lambda_1$ elements from $A_i$. It follows that to get a non-zero value, we must substitute at least $|\lambda|-\lambda_1$ elements from $J$, but $|\lambda|-\lambda_1 \geq q$, and we reach a contradiction since $J^q=\{0\}$.
	
	Suppose now that $\chi_n^L(A)=\sum\limits_{\substack{\lambda\vdash n \\ |\lambda|-\lambda_1< q}} m_\lambda \chi_\lambda$. Since $|\lambda|-\lambda_1< q$, then $\lambda_1 > n -q$ and by the hook formula $d_\lambda =\deg \chi_\lambda= \dfrac{n!}{(n-q)!}\leq n^q$. Thus by \cite[Theorem 5]{Gordienko2013JA}, it follows that
	$$c_n ^L(A)=\sum_{\substack{\lambda\vdash n \\ |\lambda|-\lambda_1< q}} m_\lambda d_\lambda \leq n^q \sum_{\substack{\lambda\vdash n \\ |\lambda|-\lambda_1< q}} m_\lambda \leq C n^{q'} $$
	for some constant $C, q'$, and the claim is proved.
\end{proof}

Next we shall give us a characterization of finite dimensional $L$-algebras with multiplicities of the $n$th differential cocharacter bounded by a constant.

We start by proving the following result.

\begin{lemma}
	\label{Lemma Mostro}
	Let $A$ be a finite dimensional $L$-algebra over an algebraically closed field such that $\dim_F A/J(A) \leq 1$. Then there exists a constant $C$ such that in $\chi_n^L(A)=\sum_{\lambda\vdash n} m_\lambda \chi_\lambda$
	$$m_\lambda \leq C,$$
	for all $n \geq 1$.
\end{lemma}

\begin{proof}
	Let $A=A_{ss}+J$ where $A_{ss}$ is a semisimple subalgebra and $J=J(A)$ is the Jacobson radical of $A$. Since $\dim_F A/J(A) \leq 1$, it follows that either $A_{ss}\cong F$ or $A = J(A)$ is a nilpotent algebra. Clearly if $A$ is a nilpotent algebra, we have nothing to prove. So let assume that $A_{ss}\cong F$. 
	
	Let now $d=\dim_F A$ and $\{a_1,\dots, a_d\}$ be a basis of $A$ where  $a_1\in A_{ss}$ and $a_2, \dots, a_d\in J$. If $q$ is the smallest positive integer such that $J^q=\{ 0\}$, we shall prove that $m_\lambda \leq d q^{dq}$ for all $\lambda \vdash n$. 
	
	 Notice that since $\dim_F A/J(A) \leq 1$, by Theorem \ref{exp}, $c_n^L(A)$ is polynomially bounded. Then, by Theorem \ref{Thm: Characterization in terms of cocharacter}, we get that $m_\lambda \neq 0$ if and only if $h(\lambda)\leq q$, where $h(\lambda)$ is the height of the partition $\lambda \vdash n$, i.e., the number of the rows of $\lambda$.  
	
	So let $\lambda \vdash n$ be a partition such that $h(\lambda)\leq q$. Consider the Young tableau $T_\lambda$ of shape $\lambda$ and the corresponding minimal essential idempotent $e_{T_\lambda}$. Then it is well-known that
	$$
	e_{T_{\lambda}}=\sum_{\substack{\sigma\in R_{T_{\lambda}} \\ \tau\in C_{T_{\lambda}}}} (\sgn\tau)\sigma\tau
	$$
	 where $R_{T_{\lambda}}$ and $C_{T_{\lambda}}$ are the subgroups of row and column permutations of $T_{\lambda}$, respectively.
	 
	 For all $1\leq j \leq q$, let $X_j$ be the set of variables whose indices lies in the $i$th row of $T_\lambda$. Thus, for any $f\in P_n^L$, the polynomial $e_{T_\lambda}f$ is symmetric in each set $X_1,\dots, X_q$ and its variables are partitioned into the disjoint union of $q$ subsets $X_1 \cup \dots \cup X_q$. Notice that $X_j$ may be empty if $h(\lambda)<j<q$. 
	 
	 Notice that for any $\rho \in S_n$, $\rho e_{T_\lambda} \neq 0$.  Then it follows that, if $e_{T_\lambda} f \neq 0$, where $f$ is a multilinear differential polynomial, then $e_{T_\lambda}f$ and $\rho e_{T_\lambda}f$ generate the same irreducible $S_n$-module.
	 
	 Let $f_1,\dots, f_m$ be a multilinear differential polynomial generating in $P_n^L(A)$ different isomorphic irreducible $S_n$-modules corresponding to the same partition. By the above, one can choose $\rho_1, \dots, \rho_m\in S_n$ and a decomposition $X=X_1\cup \dots \cup X_q$ such that $\rho_1 f_1, \dots, \rho_m f_m$ are simultaneously symmetric on $X_j$, $1\leq j \leq q$. Thus without loss of generality, we may assume that $f_1, \dots, f_m$ satisfy this condition.
	 
	 Now assume by contradiction that $m=m_\lambda >C=dq^{dq}$ and prove that $A$ satisfies a differential identity of the type
	 \begin{equation} \label{f}
	 f=\beta_1 f_1 + \dots + \beta_m f_m,
	 \end{equation}
	 where $\beta_1, \dots, \beta_m\in F$ are not all zero. Then we shall reach a contradiction since this will say that $f_1, \dots, f_m$ are linearly dependent modulo $\I^L(A)$.
	 
	 Since $f$ is multilinear, in order to verify that $f\equiv 0$, it is sufficient to verify that $f$ has only zero value on elements of a basis of $A$. First let us define substitutions of special kind. Consider the non-negative integers $\alpha_1^j, \dots, \alpha_d^j$ such that, for all $1\leq j \leq q$,
	 $$
	 \sum_{i=1}^{d}\alpha_i^j=|X_j|.
	 $$
	 We say that an evaluation $\varphi$ has type
	 $$
	 (\alpha_1^j, \dots, \alpha_d^j)
	 $$
	 for $1\leq j \leq q$, if we replace the variables from $X$ in the following way: for fixed $j$, $1\leq j \leq q$, we evaluate the first $\alpha_1^j$ variables from $X_j$ by elements $a_1$, the next $\alpha_2^j$ in $a_2$, and so on up to the last $\alpha_d^j$ variables from $X_j$ in $a_d$.
	 
	 In order to get a non-zero value of $f$ in \eqref{f}, any substitution should satisfy the following condition
	 $$
	 \sum_{i=2}^{d}\alpha_i^j\leq q-1
	 $$
	 for all $1\leq j \leq q$, since $J^q=\{ 0 \}$. Moreover, by definition we have also the following restriction
	 $$
	 \alpha_1^j=|X_j|-\sum_{i=2}^{d}\alpha_i^j
	 $$
	 for all $1\leq j \leq q$. Then for any $1\leq j \leq q$, the number of distinct $d$-tuples $(\alpha_1^j, \dots, \alpha_d^j)$ is less than $q^d$. Thus it follows that the total number $N$ of distinct type of special substitutions is less than $q^{dq}$.
	 
	 Let us consider all these $N$ distinct special substitutions $\varphi_1,\dots, \varphi_N$ and construct the matrix $(b_{ij})$, where, for all $1\leq i \leq m$ and $1\leq j \leq N$, 
	 $$
	 \varphi_j(f_i)=b_{ij}.
	 $$
	 This matrix has $m$ rows and $N$ columns of elements of $A$. Since $m> d q ^{dq} >dN$, the rows of $(b_{ij})$ are linearly dependent. Thus there exist $\beta_1, \dots \beta_m\in F$ not all zero such that
	 
	 $$
	 \sum_{i=1}^m \beta_i b_{ij}=0
	 $$
	 for all $1\leq j \leq N$, i.e., the polynomial $f=\sum_{i=1}^m \beta_i f_i$ is zero under all special substitution $\varphi_1, \dots, \varphi_N$. Therefore it is enough to show that this implies that $f\in \I^L(A)$.
	 
	 To this end, let $\psi$ be any substitution by elements of the basis $\{a_1, \dots, a_d\}$. Let $l_1^j$ be the number of variables in $X_j$ mapped by $\psi$ in $a_1$; let $l_2^j$ the number of variables in $X_j$ mapped by $\psi$ in $a_2$, and so on. Since $f$ is simultaneously symmetric on $X_1,\dots, X_q$, we get that, for all $\rho \in S_n$ such that $\rho(X_1)=X_1, \dots \rho(X_q)=X_q$, 
	 $$
	 \psi(f)=\psi(\rho f)=(\psi \rho)f.
	 $$
	 In particular, we can choose $\rho \in S_n$ such that $\psi \rho$ is the special substitution of the type $(l_1^j, \dots, l_d^j)$. By the above, $\psi(f)=(\psi \rho)f=0$ and $f\in \I^L(A)$, a contradiction. This complete the proof.
 \end{proof}

In case $A$ is a finite dimensional associative algebra we have the following result (see \cite[Sections 7.2 and 7.4]{GiambrunoZaicev2005book}).

\begin{theorem}
	\label{Theorem ordinary moultiplicities bounded}
	Let $A$ be a finite dimensional algebra over $F$. Then the following conditions are equivalent:
	\begin{itemize}
		\item[1.] $c_n(A)$ is polynomially bounded;
		\vspace{1mm}
		\item[2.] There exists a constant $C$ such that in $\chi_n(A)=\sum_{\lambda\vdash n} \bar{m}_\lambda \chi_\lambda$
		$$\bar{m}_\lambda \leq C,$$
		for all $n \geq 1$;	
		\vspace{1mm}
		\item[3.] there exists a constant $k$ such that $l_n(A)=\sum_{\lambda\vdash n} \bar{m}_\lambda \leq k$, for all $n\geq 1$.
	\end{itemize}
\end{theorem}

We are now in the position to prove the following result.

\begin{theorem}
	\label{Thm:Mult bounded by const}
	Let $L$ be a Lie algebra over a field $F$ of characteristic zero, $A$ be a finite dimensional $L$-algebra over $F$ and $\chi_n^L(A)=\sum_{\lambda\vdash n} m_\lambda \chi_\lambda$
	be its $n$th differential cocharacter. Then $c_n ^L (A)$ is polynomially bounded if and only if there exists a constant $C$ such that, for all $\lambda\vdash n$, the inequality
	$$m_\lambda \leq C$$
	holds.
\end{theorem} 

\begin{proof}
	Since the decomposition of $\chi_n^L(A)$ into irreducible characters do not change by extending the base field, we may assume that $F$ is algebraically closed. Suppose now that $c_n ^L (A)$, $n\geq 1$, is polynomially bounded, then the proof follows by Theorem \ref{Thm: Characterization in terms of equvalent algebra}, Remark \ref{Remark cocharacter and colegth} and Lemma \ref{Lemma Mostro}.
	
	Conversely, let $\chi_n^L(A)=\sum_{\lambda\vdash n} m_\lambda \chi_\lambda$
	be the $n$th differential cocharacter of $A$ and assume that there exists a constant $C$ such that, for all $\lambda\vdash n$, the inequality $m_\lambda \leq C$ holds. Then if $\chi_n(A)=\sum_{\lambda\vdash n} \bar{m}_\lambda \chi_\lambda$ is the $n$th (ordinary) cocharater of $A$, by Remark \ref{Rmk Codimensions and multiplicity} we have that $\bar{m}_\lambda \leq C$ for all $\lambda\vdash n$. Thus by Theorems \ref{Thm: Characterization in terms of ordinary exp} and \ref{Theorem ordinary moultiplicities bounded} we are done.
\end{proof}

As an important consequence, we shall prove the following corollary that relates the growth of the differential codimension sequence of a finite dimensional $L$-algebra $A$ with its differential colength.

\begin{corollary}
	\label{Cor: colength}
	Let $L$ be a Lie algebra over a field $F$ of characteristic zero and let $A$ be a finite dimensional $L$-algebra over $F$. Then $c_n ^L (A)$, $n\geq 1$, is polynomially bounded if and only if $l_n^L(A) \leq k$, for some constant $k$ and for all $n\geq 1$.
\end{corollary}

\begin{proof}
	Assume first that $c_n^L(A)$, $n\geq 1$, is polynomially bounded. By the previous theorem all non-zero multiplicities $m_\lambda$ in
	$$\chi_n^L(A)=\sum_{\lambda\vdash n} m_\lambda \chi_\lambda$$
	are bounded by a constant $C$. On the other hand, by Theorem \ref{Thm: Characterization in terms of cocharacter}, $n-\lambda_1 \leq q$ as soon as $m_\lambda \neq 0$, where $q$ is such that  $J(A)^q=\{ 0 \}$. Since the number of partition $n-\lambda_1 \leq q$ is less than $q^2$, we get
	$$l_n^L(A)=\sum_{\lambda\vdash n} m_\lambda \leq C\cdot q^2 = \mbox{const.}$$
	
	Conversely, suppose that $l_n^L(A)$ is bounded by a constant. If $\chi_n(A)=\sum_{\lambda\vdash n} \bar{m}_\lambda \chi_\lambda$ is the $n$th (ordinary) cocharacter of $A$, then as a consequence of Remark \ref{Rmk Codimensions and multiplicity} and Theorem \ref{Theorem ordinary moultiplicities bounded} we have that $l_n(A)=\sum_{\lambda\vdash n} \bar{m}_\lambda $ is bounded by a constant. Thus by theorems \ref{Thm: Characterization in terms of ordinary exp} and \ref{Theorem ordinary moultiplicities bounded}, $c_n^L (A)$ must be polynomially bounded.
\end{proof}

We now collect the results obtained in the following theorem  which gives a complete characterization of the $L$-variety generated by a finite dimensional algebras with derivations of polynomial growth.

\begin{theorem}
	Let $L$ be a Lie algebra over a field $F$ of characteristic zero and let $A$ be a finite dimensional $L$-algebra over $F$. Then the following conditions are equivalent:
	\begin{itemize}
		\item[1.] $c_n^L(A)\leq \alpha n^t$, for some constant $\alpha,t$, for all $n\geq 1$;
		
		\vspace{1mm}
		\item[2.] $\PIexp^L(A)\leq 1$;
		
			\vspace{1mm}
		\item[3.] $c_n(A)\leq \alpha n^t$, for some constant $\alpha,t$, for all $n\geq 1$;
		\vspace{1mm}
		\item[4.] $\PIexp(A)\leq 1$;
		
			\vspace{1mm}
		\item[5.]$A\sim_{T_L} B_1 \oplus \dots \oplus B_m$, with $B_1,\dots, B_m$ finite dimensional $L$-algebras over $F$ such that $\dim B_i/J(B_i)\leq 1$ for all $1\leq i \leq m$;
		
			\vspace{1mm}
		\item[6.] There exists a constant $q$ such that
		$$\chi_n^L(A)=\sum_{\substack{\lambda\vdash n \\ |\lambda|-\lambda_1< q}} m_\lambda \chi_\lambda$$
		and $J(A)^q=0$;
		
			\vspace{1mm}
		\item[7.] There exists a constant $C$ such that in $\chi_n^L(A)=\sum_{\lambda\vdash n} m_\lambda \chi_\lambda$
		$$m_\lambda \leq C$$
		for all $n \geq 1$;
		
			\vspace{1mm}
		\item[8.] there exists a constant $k$ such that $l_n^L(A)=\sum_{\lambda\vdash n} m_\lambda \leq k$ for all $n\geq 1$.
	\end{itemize}
\end{theorem}

\section*{Declarations}
Founding:
\begin{itemize}
	\item[-] C.Rizzo was supported by the Centre for Mathematics of the University of Coimbra - UIDB/00324/2020, funded by the Portuguese Government through FCT/MCTES.
	\vspace{1mm}
	\item[-] A.C. Vieira was partially supported by CNPq.
\end{itemize}
\vspace{2mm}

Competing interests: the authors have no relevant financial or non-financial interests to disclose.

\vspace{3mm}
Data transparency: data sharing not applicable to this article as no datasets were generated or analysed during the current study.

\end{document}